\documentclass[12pt]{article}
\usepackage{
latexsym, amssymb, amscd, amsthm, graphicx, 
}

 \title{{\bf Meromorphic open-string vertex algebras and Riemannian manifolds}}
 \author{Yi-Zhi Huang}
    \date{}
    \begin{document}
    \bibliographystyle{alpha}
    \maketitle
\newtheorem{thm}{Theorem}[section]
\newtheorem{defn}[thm]{Definition}
\newtheorem{prop}[thm]{Proposition}
\newtheorem{cor}[thm]{Corollary}
\newtheorem{lemma}[thm]{Lemma}
\newtheorem{rema}[thm]{Remark}
\newtheorem{app}[thm]{Application}
\newtheorem{prob}[thm]{Problem}
\newtheorem{conv}[thm]{Convention}
\newtheorem{conj}[thm]{Conjecture}
\newtheorem{cond}[thm]{Condition}
    \newtheorem{exam}[thm]{Example}
\newtheorem{assum}[thm]{Assumption}
     \newtheorem{nota}[thm]{Notation}
\newcommand{\halmos}{\rule{1ex}{1.4ex}}
\newcommand{\pfbox}{\hspace*{\fill}\mbox{$\halmos$}}

\newcommand{\nn}{\nonumber \\}

 \newcommand{\res}{\mbox{\rm Res}}
 \newcommand{\ord}{\mbox{\rm ord}}
\renewcommand{\hom}{\mbox{\rm Hom}}
\newcommand{\edo}{\mbox{\rm End}\ }
 \newcommand{\pf}{{\it Proof.}\hspace{2ex}}
 \newcommand{\epf}{\hspace*{\fill}\mbox{$\halmos$}}
 \newcommand{\epfv}{\hspace*{\fill}\mbox{$\halmos$}\vspace{1em}}
 \newcommand{\epfe}{\hspace{2em}\halmos}
\newcommand{\nord}{\mbox{\scriptsize ${\circ\atop\circ}$}}
\newcommand{\wt}{\mbox{\rm wt}\ }
\newcommand{\swt}{\mbox{\rm {\scriptsize wt}}\ }
\newcommand{\lwt}{\mbox{\rm wt}^{L}\;}
\newcommand{\rwt}{\mbox{\rm wt}^{R}\;}
\newcommand{\slwt}{\mbox{\rm {\scriptsize wt}}^{L}\,}
\newcommand{\srwt}{\mbox{\rm {\scriptsize wt}}^{R}\,}
\newcommand{\clr}{\mbox{\rm clr}\ }
\newcommand{\tr}{\mbox{\rm Tr}}
\newcommand{\C}{\mathbb{C}}
\newcommand{\Z}{\mathbb{Z}}
\newcommand{\R}{\mathbb{R}}
\newcommand{\Q}{\mathbb{Q}}
\newcommand{\N}{\mathbb{N}}
\newcommand{\CN}{\mathcal{N}}
\newcommand{\F}{\mathcal{F}}
\newcommand{\I}{\mathcal{I}}
\newcommand{\V}{\mathcal{V}}
\newcommand{\one}{\mathbf{1}}
\newcommand{\BY}{\mathbb{Y}}
\newcommand{\ds}{\displaystyle}

        \newcommand{\ba}{\begin{array}}
        \newcommand{\ea}{\end{array}}
        \newcommand{\be}{\begin{equation}}
        \newcommand{\ee}{\end{equation}}
        \newcommand{\bea}{\begin{eqnarray}}
        \newcommand{\eea}{\end{eqnarray}}
         \newcommand{\lbar}{\bigg\vert}
        \newcommand{\p}{\partial}
        \newcommand{\dps}{\displaystyle}
        \newcommand{\bra}{\langle}
        \newcommand{\ket}{\rangle}

        \newcommand{\ob}{{\rm ob}\,}
        \renewcommand{\hom}{{\rm Hom}}

\newcommand{\A}{\mathcal{A}}
\newcommand{\Y}{\mathcal{Y}}

\newcommand{\dlt}[3]{#1 ^{-1}\delta \bigg( \frac{#2 #3 }{#1 }\bigg) }

\newcommand{\dlti}[3]{#1 \delta \bigg( \frac{#2 #3 }{#1 ^{-1}}\bigg) }

\vspace{2em}



\renewcommand{\theequation}{\thesection.\arabic{equation}}
\renewcommand{\thethm}{\thesection.\arabic{thm}}
\setcounter{equation}{0} \setcounter{thm}{0} 
\maketitle

\begin{abstract}
Let $M$ be a Riemannian manifold. For $p\in M$,
the tensor algebra  of the 
negative part of the (complex) affinization 
of the tangent space  of $M$ at $p$ has a natural structure of a 
meromorphic open-string vertex algebra. These meromorphic open-string vertex algebras 
form a vector bundle over $M$ with a connection. We construct 
a sheaf $\mathcal{V}$
of meromorphic open-string vertex algebras on the sheaf of 
parallel sections of this vector bundle. Using covariant derivatives, we construct
representations on the spaces of complex smooth functions 
of the algebras of parallel tensor fields. 
These representations are used to construct a sheaf $\mathcal{W}$ of 
left $\mathcal{V}$-modules from
the sheaf of smooth functions. In particular, 
we obtain a meromorphic open-string vertex algebra $V_{M}$ of the global sections 
on $M$ of the sheaf $\mathcal{V}$ and a left $V_{M}$-module $W_{M}$ of the
global sections on $M$ of the sheaf
$\mathcal{W}$. By the definitions of meromorphic open-string vertex algebra
and left module, we obtain, among many other properties, operator product 
expansion for vertex operators. 
We also show that the Laplacian on $M$ is in fact a component of 
a vertex operator for the left $V_{M}$-module $W_{M}$ restricted to the space of 
smooth functions. 
\end{abstract}

\tableofcontents

\section{Introduction}

Conjectures by physicists on
nonlinear sigma models, especially supersymmetric nonlinear sigma 
models with Calabi-Yau manifolds as targets,  are one of the most
influential sources of inspirations and motivations
for many works in geometry in the past two or three decades. 
Classically, a nonlinear sigma model is given by the set of all 
harmonic maps from a two-dimensional Riemannian manifold to 
a Riemannian manifold (the target). The main challenge for mathematicians 
is the construction of the corresponding quantum nonlinear sigma model.
The difficulties lie in the fact that because
the target is not flat,
the nonlinear sigma model is a quantum field theory with interaction.
In physics, a quantum field theory with interaction is studied by using the 
methods of path integrals, perturbative expansion (more precisely, asymptotic 
expansion) and renormalization. Unfortunately,
it does not seem to be mathematically
possible to directly rigorize these physical methods
to construct the correlation functions for
such a quantum field theory. 

Assuming the existence of nonlinear sigma models, physicists have obtained 
many surprising mathematical conjectures. Some of these conjectures
have been proved by mathematicians using methods developed in mathematics.
But there are still many deep conjectures to be understood and 
proved. Besides proving these conjectures from physics, it is also of great importance
to understand mathematically what is going on underlying these deep conjectures. 
A mathematical construction of nonlinear sigma models would allow us to 
obtain such a deep conceptual understanding and at the same time to prove 
these conjectures. 

In the present paper, we construct meromorphic open-string vertex algebras and
their representations (see \cite{H3} for definitions and constructions)
from a Riemannian manifold. We hope that these algebras and representations will
provide a starting point for a new mathematical approach to the construction 
of nonlinear sigma models.
In the case that the target is a Euclidean space or a torus, the nonlinear sigma model
becomes a linear sigma model and can be constructed mathematically using the representations
of Heisenberg algebras. In these constructions, a crucial ingredient is
the modules for the Heisenberg algebras generated by eigenfunctions of the Laplacian
of the target. The role of the eigenfunctions can be conceptually understood
as follows: Sigma models describe perturbative string theory. When the 
strings are degenerated to points in the space, string theory becomes 
quantum mechanics. In particular, all the states in quantum mechanics 
should also be states in sigma models. 
Mathematically, quantum mechanics on a Riemannian manifold $M$ 
(without additional potential terms describing interactions) is essentially the study of 
the Schr\"{o}dinger equation 
$$i\hbar \partial_{t}\psi=\Delta \psi,$$
where $\psi$ is a function on $M\times \R$. 
Using the method of separation of variables, we first study 
a product solution $fT$ of the equation above 
where $f$ is a function  on $M$ and $T$ a function  on $\R$. Then 
there exists $\lambda\in \C$ such that 
$f$ is an eigenfunction of the Laplacian $\Delta$ with the eigenvalue $\lambda$ 
and 
$T=Ce^{-\frac{\lambda}{i\hbar}t}$
for some $C\in \C$. 
Thus the study of the Schr\"{o}dinger equation above is reduced to the study 
of eigenvalues and  eigenfunctions of the Laplacian $\Delta$. Eigenfunctions  
of $\Delta$ are states in the quantum mechanics on $M$ whose eigenvalues 
are the energies when the quantum mechanical particle is in these states.

For a Riemannian manifold $M$, its tangent spaces are Euclidean spaces. From 
these tangent spaces, one can construct vertex operator algebras associated 
with Heisenberg algebras. These vertex operator algebras form a vector bundle
of vertex operator algebras on $M$.  By tautology, the space of smooth sections of this 
bundle is a vertex algebra, a variant of vertex operator algebras
satisfying less conditions. Geometrically this vertex algebra is not 
very interesting because
as a module for the ring of the smooth functions on $M$,
the information about this vertex algebra can all be obtained from the 
theory of vector bundles and the vertex operator algebras over the fibers. 
Algebraically, since this vertex algebra does not satisfy the important grading 
restriction condition and its weight $0$ subspace is not one dimensional 
(in fact, is the infinite-dimensional space of all smooth functions), 
not many interesting results for this vertex algebra can be expected.
To obtain a vertex algebra having better properties, it is natural to consider
the subspace of parallel sections 
of this vector bundle. It was first observed by Tamanoi \cite{T1} \cite{T}
that the space of parallel sections of a vector bundle 
of vertex operator superalgebras
constructed from suitable modules for Clifford algebras has a natural structure 
of a vertex operator superalgebra. The same observation can  be made to
see the existence of a natural structure
of a vertex operator algebra on the space of parallel sections of the vector bundle
of Heisenberg vertex operator algebras mentioned above. 
However, the only functions on $M$ belonging to 
this vertex operator algebra are constant functions and, in particular, eigenfunctions
on $M$ are not in this vertex operator algebra. In fact, we do not expect that 
eigenfunctions will in general be in any vertex operator algebra because their 
eigenvalues in general are not integers. 

On the other hand, it is known that the state space of a chiral rational 
conformal field theory is mathematically the direct sum of irreducible modules for the 
chiral algebra (the vertex operator algebra of meromorphic fields) 
of the conformal field theory (see \cite{H1} and \cite{H2}). Though the nonlinear
sigma model with target $M$ is in general not even a conformal field theory, 
it would still be natural to 
look for some modules or generalized modules
that contains eigenfunctions on $M$.
To find such modules or generalized modules, 
one would have to construct a representation of the symmetric algebra on 
the tangent space at a point $p\in M$
on the space of smooth functions on an open neighborhood of $p$. 
When $M$ is not flat, however, such a representation does not exist for 
obvious reasons: If we choose a coordinate patch near $p$ and use the derivatives
with respect to the coordinates to give the representation, the representation images
of higher derivatives depends on the coordinate patch and thus
are not covariant. If we use the covariant derivatives, 
then we do not have a representation of the symmetric algebra on 
the tangent space at a point $p$;
the failure of being a representation is measured exactly by the curvature tensor. 
This failure indicates that we should consider tensor algebras 
instead of symmetric algebras.

In \cite{H3}, the author introduced a notion of meromorphic open-string vertex algebra.
A meromorphic open-string vertex algebra is an open-string vertex 
algebra in the sense of Kong and the author \cite{HK} satisfying additional 
rationality (or meromorphicity) conditions for vertex operators.  The vertex operator 
map for a meromorphic open-string vertex algebra in general does not satisfy 
the Jacobi identity, commutativity, the commutator formula, the skew-symmetry
or even the associator formula but still satisfies
rationality and 
associativity. In particular, the operator product expansion holds for vertex operators
for a meromorphic open-string vertex algebra.
In \cite{H3}, the author constructed such algebras on the tensor algebra of 
the negative part of the affinization of a vector space and left modules for 
these algebras.

In the present paper, using covariant derivatives, parallel tensor fields and 
the constructions in \cite{H3}, we construct a sheaf of  
meromorphic open-string vertex algebras from a Riemannian manifold $M$ and 
a sheaf of left modules for this sheaf generated by the space of smooth functions 
on $M$. 

More precisely, 
For a Riemannian manifold $M$, let $TM$ be the tangent bundle
of $M$, $T(TM)$ the vector bundle of the tensor algebras
on the tangent spaces at points on $M$ and $T(\widehat{TM}_{-})$ the vector bundle over $M$ 
whose fibers are the negative parts of the affinization of 
the tangent spaces of $M$. Using the meromorphic 
open-string vertex algebras constructed in \cite{H3}, 
we construct a sheaf 
$\mathcal{V}$ of  
meromorphic open-string vertex algebras on the sheaf of 
spaces of parallel sections of the vector bundle $T(\widehat{TM}_{-})$.
In particular, the space $V_{M}$ of the global sections of $\mathcal{V}$ gives a
meromorphic 
open-string vertex algebra canonically associated to $M$.
For an open subset $U$ of $M$, let $C^{\infty}(U)$ be the space of 
smooth functions on $M$. For each open subset $U$ of $M$,
we construct a representation on the space of smooth functions on $U$ 
of the algebra of parallel
sections of $T(TM)$ on $U$. Using these representations and the constructions 
of left modules for meromorphic open-string vertex algebras in \cite{H3},  we
construct a sheaf $\mathcal{W}$ of 
left modules for $\mathcal{V}$ generated by $C^{\infty}(U)$. 
In particular, the space $W_{M}$ of the global sections of $\mathcal{W}$ gives a
left $V_{M}$-module canonically associated to $M$.
By the definitions of meromorphic open-string vertex algebra
and left module, we obtain, among many other properties, operator product 
expansion for vertex operators for $V_{M}$ and $W_{M}$. 
As an example, we show that the Laplacian on $M$ is in fact a component of 
a vertex operator for the left $V_{M}$-module $W_{M}$ 
restricted to the space of smooth functions. 

The construction in the present paper can be generalized to give constructions 
of left modules generated by forms on $M$
for suitable meromorphic open-string vertex algebra associated to a Riemannian manifold $M$. 
In the case that $M$ is K\"{a}hler or Calabi-Yau, we have
stronger results. These will be discussed in future publications.

The author studied differential geometry under the supervision of
Professor Hu Hesheng as a Master student from 1982 to 1984 in Fudan university. 
The publication of the present 
paper in this issue is dedicated to the memory of 
her.  This paper was finished in 2012 and was posted 
to the arXiv on May 14, 2012. The present version
is identical to the original version except that some typos are corrected and 
two paragraphs (including this one) 
are added.  In \cite{Q1}, Qi gave the explicit examples of meromorphic 
open-string vertex
algebras and their modules associated to two-dimensional 
orientable space forms. In \cite{Q2}, to understand
modules for meromorphic open-string vertex algebras generated 
by eigenfunctions of the Laplacians on space forms, Qi obtained results
and formulated a conjecture on
covariant derivatives of such eigenfunctions. 
Research projects based on this paper
have also been actively carried out by several people including the author.  

Here the author would also like to address one issue on which 
some mathematicians and the author have different opinions. 
One opinion is that this paper is based on the parallel 
sections of vector bundles and thus cannot lead to a construction of 
the two-dimensional quantum field theory associated to a Riemannian manifold. 
For example, this opinion states that 
the conformal field theories associated to 
tori cannot be constructed based on the approach developed in this paper. 
People with this opinion obviously did not read the present paper carefully. 
In the discussion in 
this introduction above, the author has indicated clearly that, though 
the meromorphic 
open-string vertex algebra associated to a Riemannian manifold 
is obtained using parallel sections, the modules
are not. Here the author would like to 
point out another related misunderstanding about the vertex-operator-algebraic 
approach to conformal field theory. Some people mistakenly 
think that a vertex operator 
algebra determines a conformal field theory completely. This is true only 
in the case of rational conformal field theories but is wrong in general. 
One class of counterexamples to this statement is the conformal field theories
associated to irrational tori.  For all irrational tori of the same dimension, 
the associated vertex operator algebras are all the same as the 
vertex operator algebra for the corresponding Euclidean space 
(the Heisenberg vertex operator algebra in this dimension) since in this case
there is no 
larger vertex operator algebra such as the lattice vertex operator algebras 
in the rational tori case. What determines a conformal field theory associated 
to a given irrational torus is the choice of a subcategory of the category 
of modules for the Heisenberg vertex operator algebra. 
The main difficulty that the author overcome in this paper is, as discussed 
above, the construction of modules generated by eigenfunctions of 
the Laplacian. This construction is not given by parallel sections 
and uses the geometry of the Riemannian 
manifold in a crucial way. Note that eigenfunctions of the Laplacian on
a Riemannian manifold contain a lot of information about 
the Riemannian manifold. Though it has been 
known for a long time that we cannot hear the shape of 
a Riemannian manifold (that is, the eigenvalues of the Laplacian cannot
determine the Riemannian manifold up to isometries), 
eigenfunctions can indeed determine 
at least a compact Riemannian manifold since every function in a suitable 
Sobolev space can be expanded as a (finite or infinite) sum of eigenfunctions. 
Also for a torus, no matter whether it is rational or irrational, it is easy
to use the construction of the present paper to 
construct the corresponding conformal field theory. 
This is in fact one of the reasons why
the author always believes that the approach developed in this paper is correct. 

In this paper, we shall fix a Riemannian manifold $M$. 
For basic material on Riemannian geometry, we refer the reader to the book
\cite{P}. For meromorphic open-string vertex algebras and left modules, see
\cite{H3}.

The present paper is organized as follows: In Section 2, we recall
some basic constructions of 
vector bundles and sheaves on a Riemannian manifold $M$. 
In Section 3, 
we construct the sheaf $\mathcal{V}$ of meromorphic open-string vertex algebras on $M$.
In particular, we construct the meromorphic open-string vertex algebra $V_{M}$ of 
the global sections  of $\mathcal{V}$ canonically associated to $M$.
In section 4, using covariant derivatives, we construct a homomorphism of 
algebras from the algebra of parallel tensor fields on an open subset of $M$ 
to the algebra of linear operators
on the space of smooth functions on the same open subset. In particular, we 
obtain a representation on the space of smooth functions of the 
algebra of parallel tensor fields. 
We construct in Section 5 the sheaf $\mathcal{W}$ of 
left modules for $\mathcal{V}$ generated by the sheaf of smooth functions
on $M$. In particular, we construct the left $V_{M}$-module 
$W_{M}$ of $\mathcal{W}$ canonically associated to $M$. In particular, 
we construct the left $V_{M}$-module of the 
global sections of  $\mathcal{W}$ canonically associated to $M$.
In Section 6, we show that the Laplacian on $M$ is in fact a component of 
a vertex operator for the left $V_{M}$-module $W_{M}$ restricted to 
the space of smooth functions.

\paragraph{Acknowledgments}
The author was supported in part by NSF
grant PHY-0901237.

\section{Vector bundles and sheaves from the tangent bundle
of a Riemannian manifold $M$}

In this section, we recall some basic constructions of 
vector bundles and sheaves on a Riemannian manifold.

In the present paper, we shall work with vector spaces 
over $\R$ and with vector spaces over $\C$. 
We shall use $\otimes_{\R}$ to denote the tensor product bifunctor
for the category of vector spaces over $\R$ but use $\otimes$
(omitting the subscript $\C$) to denote the tensor product bifunctor
for the category of vector spaces over $\C$. We shall use the same 
notations to denote tensor products of vector bundles 
with vector spaces over $\R$ and $\C$ as fibers. 
For a vector space $V$ over $\R$, we shall use 
$V^{\C}$ to denote the complexification $V\otimes_{\R}\C$.
For a vector bundle $E$ with vector spaces over $\R$ as fibers, 
we use $E^{\C}$ to denote the vector bundle obtained from 
$E$ by complexifying its fibers.

Let $M$ be a Riemannian manifold and $g$ the metric
on $M$. Consider the tangent bundle $TM$ of $M$ and
the trivial bundles $M\times \C[t, t^{-1}]$ and $M\times \C \mathbf{k}$ where 
$t$ is a formal variable and $\mathbf{k}$ is a basis of a one-dimensional 
vector space $\C \mathbf{k}$.                                                  
Let 
$$\widehat{TM}=(TM\otimes_{\R} (M\times \C[t, t^{-1}]))\oplus M\times \C \mathbf{k}.$$ 
be the vector bundle whose fiber at $p\in M$ is
$$\widehat{T_{p}M}=(T_{p}M\otimes_{\R} \C[t, t^{-1}])\oplus \C \mathbf{k}.$$
Since $\widehat{T_{p}M}$ for $p\in M$ has a structure of Heisenberg algebra and 
the transition functions at points of $M$ 
preserve the gradings of the Heisenberg algebras,
$\widehat{TM}$ has a structure of a vector bundle of Heisenberg algebras.
For $p\in M$, $\widehat{T_{p}M}$ has a decomposition 
$$\widehat{T_{p}M}= \widehat{T_{p}M}_{-}\oplus \widehat{T_{p}M}_{0}
\oplus \widehat{T_{p}M}_{+},$$
where 
\begin{eqnarray*}
\widehat{T_{p}M}_{-}&=&T_{p}M\otimes_{\R} t^{-1}\C[t^{-1}],\\
\widehat{T_{p}M}_{0}&=&(T_{p}M\otimes_{\R} \C t^{0})\oplus \C \mathbf{k}\nn
&\simeq& T_{p}M^{\C}\oplus \C \mathbf{k},\\
\widehat{T_{p}M}_{+}&=&T_{p}M\otimes_{\R} t\C[t].
\end{eqnarray*}
These triangle decompositions of the Heisenberg algebras give the triangle 
decomposition 
$$\widehat{TM}= \widehat{TM}_{-}\oplus \widehat{TM}_{0}\oplus \widehat{TM}_{+},$$
where 
\begin{eqnarray*}
\widehat{TM}_{-}&=&TM\otimes_{\R} (M\times t^{-1}\C[t^{-1}]),\\
\widehat{TM}_{0}&=&(TM\otimes_{\R} (M\times \C t^{0}))
\oplus M\times \C \mathbf{k}\nn
&\simeq& TM^{\C}\oplus (M\times \C \mathbf{k}),\\
\widehat{TM}_{+}&=&TM\otimes_{\R} (M\times t\C[t]).
\end{eqnarray*}

The connection on $TM$ induces connections on $\widehat{TM}$, $\widehat{TM}_{-}$ and 
$\widehat{TM}_{+}$. The product bundle $M\times \C \mathbf{k}$ has a 
trivial connection.

For $p\in M$, recall the subalgebra $N(\widehat{T_{p}M})$ of the tensor algebra
$T(\widehat{T_{p}M})$ introduced in Section 3 of \cite{H3}. In fact, let 
$I$ be the two-sided ideal of $T(\widehat{T_{p}M})$ generated by elements of the form
\begin{eqnarray*}
&(X\otimes_{\R} t^{m})\otimes (Y\otimes_{\R} t^{n})
- (Y\otimes_{\R} t^{n})\otimes (X\otimes_{\R} t^{m})
-m(a, b)\delta_{m+n, 0}\mathbf{k},&\\
&(X\otimes_{\R} t^{k})\otimes (Y\otimes_{\R} t^{0})
-(Y\otimes_{\R} t^{0})\otimes (X\otimes_{\R} t^{k}),&\\
&(X\otimes_{\R} t^{k})\otimes \mathbf{k}-\mathbf{k}\otimes (X\otimes_{\R} t^{k})&
\end{eqnarray*}
for $X, Y\in T_{p}M$, $m\in \Z_{+}$, $n\in -\Z_{+}$, $k\in \Z$. 
Then by Proposition 3.1 in \cite{H3}, 
$$N(\widehat{T_{p}M})=T(\widehat{T_{p}M})/I$$
is isomorphic to 
\begin{equation}\label{PBW}
T(\widehat{T_{p}M}_{-})\otimes T(\widehat{T_{p}M}_{+}) \otimes T(T_{p}M^{\C})
\otimes 
T(\C\mathbf{k}).
\end{equation} 
Let 
$$T(\widehat{TM}_{-}),\; 
T(\widehat{TM}_{+}), \; T(TM^{\C}), \;T(M\times \C\mathbf{k})$$ 
be the vector bundles whose fibers at $p\in M$ are
the tensor algebras 
$$T(\widehat{T_{p}M}_{-}),\; 
T(\widehat{T_{p}M}_{+}), \; T(T_{p}M^{\C}), \; T(\C\mathbf{k})$$
on the fibers of 
$$\widehat{TM}_{-},\; \widehat{TM}_{+},\;
TM^{\C}, \;M\times \C
\mathbf{k},$$ 
respectively.
Since (\ref{PBW}) is the fiber of the vector bundle 
\begin{equation}\label{PBW-bundle}
T(\widehat{TM}_{-})\otimes T(\widehat{TM}_{+}) \otimes T(TM^{\C})\otimes 
T(M\times\C\mathbf{k})
\end{equation} 
at $p\in M$, we also have a vector bundle 
$N(\widehat{TM})$ whose fiber at $p\in M$ is $N(\widehat{T_{p}M})$. 
By definition, $N(\widehat{TM})$ as a vector bundle
is isomorphic to the vector bundle (\ref{PBW-bundle}).

For a vector bundle $E$ over $M$, we shall use $\Gamma_{U}(E)$ 
to denote the space of smooth sections of $E$ on an open subset $U$ of $M$. For 
a vector bundle $E$ over $M$ with a connection, we shall use $\Pi_{U}(E)$
to denote the space of parallel sections of $E$ on $U$.
By definition, 
$\Pi_{U}(E)\subset \Gamma_{U}(E)$. 
When the fibers of $E$ are 
associative algebras, $\Gamma_{U}(E)$ has a structure of an associative algebra. 
If the covariant derivative with respect every element of $\Gamma_{U}(TM^{\C})$ 
is a derivation of the associative algebra $\Gamma_{U}(E)$,
then $\Pi_{U}(E)$ is  a subalgebra of $\Gamma_{U}(E)$.

Taking $E$ to be 
\begin{equation}\label{v-bundles}
N(\widehat{TM}),\; 
T(\widehat{TM}_{-}), \; T(\widehat{TM}_{+}),\;T(TM^{\C}),\;
T(M\times  \C\mathbf{k}),
\end{equation}
we have the associative algebras
\begin{eqnarray}\label{sections}
&\Gamma_{U}(N(\widehat{TM})), \;
\Gamma_{U}(T(\widehat{TM}_{-})),\; \Gamma_{U}(T(\widehat{TM}_{+})),&\nn
&\Gamma_{U}(T(TM^{\C})), \;\Gamma_{U}(T(M\times \C \mathbf{k})), &
\end{eqnarray}
respectively, of smooth sections. 
It is clear that
$$\Gamma_{U}(T(\widehat{TM}_{-})),
\Gamma_{U}(T(\widehat{TM}_{+})),
\Gamma_{U}(T(TM^{\C})), \;
\Gamma_{U}(T(M\times \C \mathbf{k}))$$
can be embedded
as subalgebras of $\Gamma_{U}(N(\widehat{TM}))$.
The connections on $\widehat{TM}_{-}$, $TM^{\C}$ and 
$\widehat{TM}_{+}$ uniquely determine connections on 
$T(\widehat{TM}_{-})$, $T(TM^{\C})$ and $T(\widehat{TM}_{+})$, respectively,
by requiring that for every open subset $U$ of $M$, 
the covariant derivatives with respect to 
every element of $\Gamma_{U}(TM)$ 
are derivations of the associative algebras
$\Gamma_{U}(T(\widehat{TM}_{-}))$, $\Gamma_{U}(T(\widehat{TM}_{+}))$ and 
$\Gamma_{U}(T(TM^{\C}))$, respectively.
We also have a canonical flat connection 
on the trivial bundle 
$$T(M\times \C \mathbf{k})\simeq M\times T(\C \mathbf{k}).$$ 
Since $N(\widehat{TM})$ is isomorphic to (\ref{PBW-bundle}),
the connections on $T(\widehat{TM}_{-})$, $T(\widehat{TM}_{+})$, $T(TM^{\C})$ and
$T(M\times \C \mathbf{k})$ further determine a connection on 
$N(\widehat{TM})$. 

By definition, the covariant derivatives 
with respect to elements of the space $\Gamma_{U}(TM)$ of the vector 
bundles in (\ref{v-bundles}) are derivations of the corresponding associative algebras 
in (\ref{sections}).
Thus we have the associative algebras
\begin{eqnarray*}
&\Pi_{U}(N(\widehat{TM})),\; 
\Pi_{U}(T(\widehat{TM}_{-})), &\nn
&
\Pi_{U}(T(\widehat{TM}_{+})),\; \Pi_{U}(T(TM^{\C})),\;
\Pi_{U}(M\times T(\C \mathbf{k}))&
\end{eqnarray*}
of parallel sections. 

For a vector bundle $E$, the spaces $\Gamma_{E}(U)$ of smooth sections on 
open subsets $U$ of $M$ and the obvious restriction maps from 
$\Gamma_{E}(U)$ to $\Gamma_{E}(U')$ when $U'\subset U$ give a sheaf
$\Gamma_{E}$. 
Similarly for a vector bundle $E$ with a connection, we also have the sheaf 
$\Pi_{E}$ whose sections on an open subset $U$ is $\Pi_{E}(U)$. 
The sheaf $\Pi_{E}$ is a subsheaf of $\Gamma_{E}$. 
Taking $E$ to be the vector bundles in (\ref{v-bundles}), we have the sheaves
\begin{eqnarray*}
&\Gamma(N(\widehat{TM})),\;
\Gamma(T(\widehat{TM}_{-})), \; \Gamma(T(\widehat{TM}_{+})),\nn
&\Gamma(T(TM^{\C})),\;
\Gamma(T(M\times \C \mathbf{k})), &\nn
&\Pi(N(\widehat{TM})),\;
\Pi(T(\widehat{TM}_{-})),\; \Pi(T(\widehat{TM}_{+})),&\nn 
&\Pi(T(TM^{\C})),\;
\Pi(M\times T(\C \mathbf{k})).&
\end{eqnarray*}

We know that the space of parallel sections of a vector bundle with a connection 
is canonically isomorphic to the space of fixed points of 
a fiber under the action of the holonomy group. In particular, we have:

\begin{prop}\label{parallel-fixed-pts}
Let $U$ be an open subset of $M$. The space 
$$\Pi_{U}(T(\widehat{TM}_{-})), \;\Pi_{U}(T(\widehat{TM}_{+})), \;
\Pi_{U}(T(TM^{\C})),\; \Pi_{U}(N(\widehat{TM}))$$ 
are
canonically isomorphic to the spaces of fixed points of 
$$T(\widehat{T_{p}M}_{-}), \; T(\widehat{T_{p}M}_{+}), \; T(T_{p}M^{\C}),\;
N(\widehat{T_{p}M}),$$ 
respectively, for $p\in U$ under the
actions of the holonomy groups of the restrictions of the vector bundles 
$$T(\widehat{TM}_{-}), \; T(\widehat{TM}_{+}), \; T(TM^{\C}),\; 
N(\widehat{TM}),$$
respectively, to $U$. \epf
\end{prop}

\section{A sheaf $\mathcal{V}$ of meromorphic open-string vertex algebras on $M$}

In this section, we construct
a sheaf of meromorphic open-string vertex algebras on $M$. In particular, the 
global 
sections of this sheaf give a canonical meromorphic open-string vertex algebra
associated to $M$. 

First we have:

\begin{prop}\label{buncle-mero-op-va}
The fibers of the vector bundle $T(\widehat{TM}_{-})$ have natural 
structures of meromorphic open-string vertex algebras and 
$T(\widehat{TM}_{-})$ has a natural structure of  
vector bundle of meromorphic open-string vertex algebras.
\end{prop}
\pf
Since the fibers of $T(\widehat{TM}_{-})$ 
are the tensor algebras on the fibers of $\widehat{TM}_{-}$,
by Theorem 5.1 in \cite{H3},  they 
have natural structures of meromorphic open-string vertex algebras. 
It is clear that the transition functions of the vector bundle $T(\widehat{TM}_{-})$
at points on $M$ are automorphisms of meromorphic open-string vertex algebras.
Thus $T(\widehat{TM}_{-})$ has a natural structure of
vector bundle of meromorphic open-string vertex algebras. 
\epfv

Then we have:

\begin{cor}\label{sec-mero-op-va}
For an open subset $U$ of $M$, 
the space $\Gamma_{U}(T(\widehat{TM}_{-}))$ of sections of $T(\widehat{TM}_{-})$ 
has a natural structure of meromorphic open-string vertex algebra. The assignment
$$U\to \Gamma_{U}(T(\widehat{TM}_{-}))$$
together with the restrictions of sections
form a sheaf of meromorphic open-string vertex algebras.
\end{cor}
\pf
The $\Z$-gradings on the fibers of $T(\widehat{TM}_{-})$ induce a $\Z$-grading on 
$\Gamma_{U}(T(\widehat{TM}_{-}))$. The constant section $1$ is the vacuum.
The vertex operator map is defined pointwise. It is clear that 
with the $\Z$-grading, the vacuum and the vertex operator map, 
$\Gamma_{U}(T(\widehat{TM}_{-}))$ is a meromorphic open-string vertex algebra. 
The second conclusion is also clear.
\epfv

The construction in 
Corollary \ref{sec-mero-op-va} is simple. But these meromorphic open-string vertex algebras
are not what we are interested in. In fact, the sheaf of 
meromorphic open-string vertex algebras
obtained in Corollary \ref{sec-mero-op-va} contains the sheaf of smooth functions
on $M$ and the smooth functions commute with vertex operators. In particular, 
the vertex operators in this sheaf of meromorphic open-string vertex algebras
cannot contain differential operators acting on the space of smooth functions. 
Since the quantum mechanics on $M$ involves differential operators, 
the sheaf of meromorphic open-string vertex algebras 
in Corollary \ref{sec-mero-op-va} is not what we are looking for. 
We shall instead construct a sheaf of meromorphic open-string vertex algebras
using parallel sections.

Given a meromorphic open-string vertex algebra $(V, Y_{V}, \one)$ 
and a group $H$ of automorphisms of 
$V$, let $V^{H}$ be the subspace of $V$ consisting of elements that are fixed by $H$.
Since automorphisms of $V$ preserve $\one\in V$ (see \cite{H3}), 
$\one \in V^{H}$. Also since for $u, v\in V^{G}$ and $h\in H$, 
$hY_{V}(u, x)v=Y_{V}(hu, x)hv=Y_{V}(u, x)v$, the image of $V^{H}\otimes V^{H}$
under $Y_{V}$ is in $V^{H}[[x, x^{-1}]]$. We shall denote the restriction of  $Y_{V}$
to $V^{H}\otimes V^{H}$ by $Y_{V^{H}}$. Then $Y_{V^{G}}$ is a linear map 
from $V^{H}\otimes V^{H}$ to $V^{H}[[x, x^{-1}]]$. The following result is 
obvious:

\begin{prop}
The triple $(V^{H}, Y_{V^{H}}, \one)$  is a meromorphic open-string vertex subalgebra
of $(V, Y_{V}, \one)$. \epf
\end{prop}

For $p\in M$ and a connected open subset $U$ of $M$ containing $p$, 
the holonomy group $H_{p}(U)$ of the 
restriction of the vector bundle $T(\widehat{TM}_{-})$ to $U$
acts on the 
fiber $T(\widehat{T_{p}M}_{-})$ at $p$ of the vector bundle $T(\widehat{TM}_{-})$. 
By Proposition \ref{buncle-mero-op-va}, 
$T(\widehat{T_{p}M}_{-})$ has a structure of meromorphic open-string vertex algebra.

\begin{lemma}
For a connected open subset $U$ of $M$, 
$\alpha\in H_{p}(U)$ and $u, v\in T(\widehat{T_{p}M}_{-})$, 
$$\alpha(Y_{T(\widehat{T_{p}M}_{-})}(u, x)v)
=Y_{T(\widehat{T_{p}M}_{-})}(\alpha(u), x)\alpha(v).$$
\end{lemma}
\pf
Recall the notations in \cite{H3}. We need only prove the lemma in the case 
$$u=X_{1}(-n_{1})\cdots X_{k}(-n_{k})\mathbf{1}$$
for $X_{1}, \dots, X_{k}\in T_{p}M$ and $n_{1}, \dots, n_{k}\in \Z_{+}$.
Since the connection on $T(\widehat{TM}_{-})$ is induced 
from the connection on $TM^{\C}$, the parallel transport in $T(\widehat{TM}_{-})$
along a path in $M$ is also induced from the parallel transport 
in $TM^{\C}$ along the same path. Let $\gamma$ be a loop in $M$ based at $p$.
Denote both the parallel transports along $\gamma$ in $TM^{\C}$ and in $T(\widehat{TM}_{-})$
by $\alpha_{\gamma}$. Then we have 
$$\alpha_{\gamma}(\one)=\one$$
and 
\begin{equation}\label{alpha-elements}
\alpha_{\gamma}(X_{1}(-m_{1})\cdots X_{k}(-m_{k})\mathbf{1})
=\alpha_{\gamma}(X_{1})(-m_{1})\cdots \alpha_{\gamma}(X_{k})(-m_{k})\mathbf{1}
\end{equation}
for $n_{1}, \dots, n_{k}\in \Z$.

By definition,
\begin{eqnarray*}
\lefteqn{Y_{T(\widehat{T_{p}M}_{-})}(X_{1}(-n_{1})\cdots X_{k}(-n_{k})\one, x)v}\nn
&&=\left(\nord \frac{1}{(n_{1}-1)!}\left(\frac{d^{n_{1}-1}}{dx^{n_{1}-1}}
X_{1}(x)\right)\cdots \frac{1}{(n_{k}-1)!}
\left(\frac{d^{n_{k}-1}}{dx^{n_{k}-1}}
X_{k}(x)\right)\nord \right) v, 
\end{eqnarray*}
where as in \cite{H3}, 
$$X_{i}(x)=\sum_{n\in \Z}X_{i}(n)x^{-n-1}$$
for $i=1, \dots, k$ and $\nord \cdot \nord$ is the normal ordering
operation defined in \cite{H3}. Thus by Lemma 4.2 and (\ref{alpha-elements}), we have 
\begin{eqnarray*}
\lefteqn{\alpha_{\gamma}(Y_{T(\widehat{T_{p}M}_{-})}(u, x)v)}\nn
&&=\alpha_{\gamma}(Y_{T(\widehat{T_{p}M}_{-})}(X_{1}(-n_{1})
\cdots X_{k}(-n_{k})\one, x)v)\nn
&&=\alpha_{\gamma}\left(\left(\nord \frac{1}{(n_{1}-1)!}\left(\frac{d^{n_{1}-1}}{dx^{n_{1}-1}}
X_{1}(x)\right)\cdots \frac{1}{(n_{k}-1)!}
\left(\frac{d^{n_{k}-1}}{dx^{n_{k}-1}}
X_{k}(x)\right)\nord \right) v\right)\nn
&&=\nord \frac{1}{(n_{1}-1)!}\left(\frac{d^{n_{1}-1}}{dx^{n_{1}-1}}
(\alpha_{\gamma}(X_{1}))(x)\right)\cdot\nn
&&\quad\quad\quad\quad\quad\quad \cdots \frac{1}{(n_{k}-1)!}
\left(\frac{d^{n_{k}-1}}{dx^{n_{k}-1}}
(\alpha_{\gamma}(X_{k}))(x)\right)\nord \alpha_{\gamma}(v)\nn
&&=Y_{T(\widehat{T_{p}M}_{-})}(\alpha_{\gamma}(X_{1})(-n_{1})\cdots 
\alpha_{\gamma}(X_{k})(-n_{k})\one, x)\alpha_{\gamma}(v)\nn
&&=Y_{T(\widehat{T_{p}M}_{-})}(\alpha_{\gamma}(X_{1}(-n_{1})\cdots 
X_{k}(-n_{k})\one), x)\alpha_{\gamma}(v)\nn
&&=Y_{T(\widehat{T_{p}M}_{-})}(\alpha_{\gamma}(u), x)\alpha_{\gamma}(v).
\end{eqnarray*}
\epfv

From the lemma above, we obtain immediately:

\begin{cor}\label{holonomy}
For a connected open subset $U$ of $M$, the holonomy group 
$H_{p}(U)$ is a subgroup of the automorphism group 
of the meromorphic open-string vertex algebra $T(\widehat{T_{p}M}_{-})$. In particular, 
$(T(\widehat{T_{p}M}_{-}))^{H_{p}(U)}$ is a meromorphic open-string vertex subalgebra
of $T(\widehat{T_{p}M}_{-})$.
\epf
\end{cor}

For an open subset $U$ of $M$, 
let 
$$V_{U}=\Pi_{U}(T(\widehat{TM}_{-})).$$
Then the assignment $U\to V_{U}$ and the 
restrictions of sections give a sheaf $\mathcal{V}$. 
By Proposition \ref{parallel-fixed-pts}, $V_{U}$
is canonically isomorphic to $(T(\widehat{T_{p}M}_{-}))^{H_{p}(U)}$. Thus we have:

\begin{thm}\label{V-U-mero-op-va}
For a connected open subset $U$ of $M$ and $p\in U$, the canonical isomorphism from
$(T(\widehat{T_{p}M}_{-}))^{H_{p}(U)}$ to 
$V_{U}$ gives $V_{U}$ a natural structure of meromorphic open-string vertex algebra. 
This structure of meromorphic open-string vertex algebra is independent of the choice of $p$.
For general open subset $U$ of $M$, $V_{U}$ as a $\Z$-graded vector space 
is isomorphic to the underlying $\Z$-graded vector space
of the direct product meromorphic open-string vertex algebra 
$\prod_{\alpha\in \mathcal{A}}V_{U_{\alpha}}$ (see Definition 2.6 
in \cite{H3}) where $U_{\alpha}$ for
$\alpha\in \mathcal{A}$ are the connected components of $U$. In particular, 
$V_{U}$ also has a natural structure of meromorphic open-string vertex algebra.
For an open subset $U$ of $M$ and an open subset $\tilde{U}$ of $U$, the restriction map 
from $V_{U}$ to $V_{\tilde{U}}$ is a homomorphism of meromorphic open-string vertex algebras.
In particular, the sheaf 
$\mathcal{V}$ is a sheaf of meromorphic open-string vertex algebras.
\end{thm}
\pf
The first and second statements of the theorem are clear. 

For general open subset $U$ of $M$,
choose a point $p_{\alpha}$ in each connected component $U_{\alpha}$ of $U$ for 
$\alpha\in \mathcal{A}$ (elements of $\mathcal{A}$ 
labeling the connected components of $U$), then
$\Pi_{U_{\alpha}}(T(\widehat{T_{p_{\alpha}}M}_{-}))$ is isomorphic to 
$(T(\widehat{T_{p_{\alpha}}M}_{-}))^{H_{p_{\alpha}}(U_{\alpha})}$
as  a graded vector space, where $H_{p_{\alpha}}(U_{\alpha})$ is 
the holonomy group of of the connection on 
the vector bundle $T(\widehat{TM}_{-})$ restricted to the 
connected component $U_{\alpha}$. But $\Pi_{U}(T(\widehat{TM}_{-}))$ 
is isomorphic to $\prod_{\alpha\in \mathcal{A}} 
\Pi_{U_{\alpha}}(T(\widehat{T_{p_{\alpha}}M}_{-}))$
as  a graded vector space.
Hence $V_{U}$ is isomorphic to 
$\prod_{\alpha\in \mathcal{A}} 
(T(\widehat{T_{p_{\alpha}}M}_{-}))^{H_{p_{\alpha}}(U_{\alpha})}$
as  a graded vector space.
Since $\prod_{\alpha\in \mathcal{A}} 
(T(\widehat{T_{p_{\alpha}}M}_{-}))^{H_{p_{\alpha}}(U_{\alpha})}$ 
has a  structure of the direct product  meromorphic open-string vertex algebra of 
$(T(\widehat{T_{p_{\alpha}}M}_{-}))^{H_{p_{\alpha}}(U_{\alpha})}$ 
for $\alpha\in \mathcal{A}$, $V_{U}$ has a natural structure of 
a meromorphic open-string vertex algebra.

For an open subset $U$ of $M$ and an open subset $\tilde{U}$ of $U$, let 
$U_{\alpha}$ for $\alpha\in \mathcal{A}$ be the connected components 
of $U$ and let $\tilde{U}_{\beta}$ for $\beta\in \mathcal{B}$
be the connected components of $\tilde{U}$. Then for $\beta\in \mathcal{B}$, 
there exists $\alpha\in \mathcal{A}$ 
such that $\tilde{U}_{\beta}\subset U_{\alpha}$. For each $\beta\in \mathcal{B}$, 
we choose a point $\tilde{p}_{\beta}\in \tilde{U}_{j}$. Then there 
exists $\alpha\in \mathcal{A}$ such that 
$\tilde{p}_{\beta}\in U_{\alpha}$. We choose 
$p_{\alpha}\in U_{\alpha}$ from those $\tilde{p}_{\beta}$'s
such that  $\tilde{p}_{\beta}\in \tilde{U}_{\beta}$. Then
$H_{\tilde{p}_{\beta}}(\tilde{U}_{\beta})$ can be naturally 
embedded into $H_{p_{\alpha}}(U_{\alpha})$
when $\tilde{p}_{\beta}\in U_{\alpha}$. 
Thus the direct product meromorphic open-string vertex algebra
$\prod_{\alpha\in \mathcal{A}} 
(T(\widehat{T_{p_{\alpha}}M}_{-}))^{H_{p_{\alpha}}(U_{\alpha})}$ can be embedded into 
the direct product meromorphic open-string vertex algebra
$\prod_{\beta\in \mathcal{B}} 
(T(\widehat{T_{\tilde{p}_{\beta}}M}_{-}))^{H_{\tilde{p}_{\beta}}(\tilde{U}_{\beta})}$. 
The embedding 
from $\prod_{\alpha\in \mathcal{A}} 
(T(\widehat{T_{p_{\alpha}}M}_{-}))^{H_{p_{\alpha}}(U_{\alpha})}$
to 
$\prod_{\beta\in \mathcal{B}} 
(T(\widehat{T_{\tilde{p}_{\beta}}M}_{-}))^{H_{\tilde{p}_{\beta}}(\tilde{U}_{\beta})}$
corresponds to 
the restriction map from $V_{U}$ to $V_{\tilde{U}}$, that is, we have the following 
commutative diagram:
$$\begin{CD}
\prod_{\alpha\in \mathcal{A}} 
(T(\widehat{T_{p_{\alpha}}M}_{-}))^{H_{p_{\alpha}}(U_{\alpha})} @>>> V_{U}\\
@VVV @VVV\\
\prod_{\beta\in \mathcal{B}} 
(T(\widehat{T_{\tilde{p}_{\beta}}M}_{-}))^{H_{\tilde{p}_{\beta}}(\tilde{U}_{\beta})}
 @>>> V_{\tilde{U}}.
\end{CD}$$
Since the embedding 
from $\prod_{\alpha\in \mathcal{A}} 
(T(\widehat{T_{p_{\alpha}}M}_{-}))^{H_{p_{\alpha}}(U_{\alpha})}$
to 
$\prod_{\beta\in \mathcal{B}} 
(T(\widehat{T_{\tilde{p}_{\beta}}M}_{-}))^{H_{\tilde{p}_{\beta}}(\tilde{U}_{\beta})}$
is a homomorphism of meromorphic open-string vertex algebras, the restriction map 
from $V_{U}$ to $V_{\tilde{U}}$ is also a homomorphism of 
meromorphic open-string vertex algebras.
\epfv

\begin{rema}
{\rm For an open subset $U$ of $M$, $V_{U}$ is always nontrivial.
In fact. the metric $g$ can be viewed as an element of the space $\Gamma_{U}(T^{2}(T^{*}M))$
of smooth sections on $U$ of the second symmetric tensor powers of the cotangent bundle
$T^{*}M$ of $M$. On the other hand, $g$ also
gives an isomorphism of vector bundles from $T^{*}M$ to $TM$.
It induces an isomorphism of vector bundles from $T^{2}(T^{*}M)$ to $T^{2}(TM)$,
which in turn induces a linear isomorphism from $\Gamma_{U}(T^{2}(T^{*}M))$
to $\Gamma_{U}(T^{2}(TM))$. The image of the element $g\in \Gamma_{U}(T^{2}(T^{*}M))$
under this isomorphism is an element of $\Gamma_{U}(T^{2}(TM))$ and 
thus gives an element of $\Gamma_{U}(T^{2}(TM^{\C}))$. We denote this element 
of $\Gamma_{U}(T^{2}(TM^{\C}))$ by
$g_{\C}^{-1}$. Since $g$ is parallel, $g_{\C}^{-1}$ is also parallel. 
For $k, l\in \Z_{+}$, the vector bundles $TM\otimes (M\times \C t^{-k})$ and 
$TM\otimes (M\times \C t^{-l})$
are isomorphic to $TM^{\C}$. In particular, the space 
$$\Gamma_{U}((TM\otimes (M\times \C t^{-k}))\otimes
(TM\otimes (M\times \C t^{-l})))$$
of sections of the vector bundle 
$$(TM\otimes (M\times \C t^{-k}))\otimes
(TM\otimes (M\times \C t^{-l}))$$
is isomorphic to the space $\Gamma_{U}(T^{2}(TM^{\C}))$. 
In particular, $g_{\C}^{-1}\in \Gamma_{U}(T^{2}(TM^{\C}))$ corresponds to an element
$$g_{\C}^{-1}(-k, -l)\in \Gamma_{U}((TM\otimes (M\times \C t^{-k}))\otimes
(TM\otimes (M\times \C t^{-l})).$$
Since $g_{\C}^{-1}$ is parallel and the connection on 
$(TM\otimes (M\times \C t^{-k}))\otimes
(TM\otimes (M\times \C t^{-l}))$ is induced from the 
connection on $T^{2}(TM^{\C})$, 
$g_{\C}^{-1}(-k, -l)$ is also parallel, that is, 
\begin{eqnarray*}
g_{\C}^{-1}(-k, -l)&\in& 
\Pi_{U}((TM\otimes (M\times \C t^{-k}))\otimes
(TM\otimes (M\times \C t^{-l})))\nn
&\subset& \Pi_{U}(T(\widehat{TM}_{-}))\nn
&=&V_{U},
\end{eqnarray*}
for $k, l\in \Z_{+}$, giving infinitely many nonzero elements of $V_{U}$ of 
different weights. }
\end{rema}

\begin{rema}
{\rm It is well known that for $p\in M$, the symmetric algebra 
$S(\widehat{T_{p}M}_{-})$ has a natural structure of a vertex operator algebra. 
These symmetric algebras form a vector bundle $S(\widehat{TM}_{-})$ 
of vertex operator algebras with a connection. The same construction 
as the one for $\mathcal{V}$ above shows that the spaces
$\Pi_{U}(S(\widehat{TM}_{-})$
of parallel sections of $S(\widehat{TM}_{-})$ on open subsets $U$ of $M$ form 
a sheaf of conformal vertex algebras such that when $U$ is connected, 
$\Pi_{U}(S(\widehat{TM}_{-})$ is a vertex operator algebra. 
From Remark 5.2 in \cite{H3}, for $p\in M$,
we have a homomorphism of meromorphic open-string vertex algebras
from $T(\widehat{T_{p}M}_{-})$ to $S(\widehat{T_{p}M}_{-})$. Thus we have 
a homomorphism of vector bundles from $T(\widehat{TM}_{-})$ to
$S(\widehat{TM}_{-})$ such that the connection on $T(\widehat{TM}_{-})$
is mapped to the connection on $S(\widehat{TM}_{-})$. In particular, we have a 
homomorphism of sheaves of meromorphic open-string vertex algebras 
from the sheaf $\mathcal{V}$ to the sheaf of the spaces
of parallel sections
of the vector bundle $S(\widehat{TM}_{-})$. }
\end{rema}

\section{Covariant derivatives and parallel tensor fields}

Given an open subset $U$ of $M$, let $C^{\infty}(U)$ be the space of 
complex smooth functions 
on $U$. For $m\in \N$, let $T^{m}(TM^{\C})$ be the $m$-th tensor power of 
the tangent bundle $TM^{\C}$ and $\Gamma_{U}(T^{m}(TM^{\C}))$ the space of sections of 
$T^{m}(TM^{\C})$. Then $\Gamma_{U}(T(TM^{\C}))$ is the coproduct of 
$\Gamma_{U}(T^{m}(TM^{\C}))$, $m\in \N$. 
Given $f\in C^{\infty}(U)$, there is an $m$-th order 
covariant derivative $\nabla^{m} f$ which can be viewed as a $(0, m)$-tensor. 
Note that $\nabla^{m} f$
is originally defined only for real smooth function $f$ and applies only to 
real tensor fields. But it can be extended in the obvious way to a
complex smooth function $f$ and applies to complex tensor fields. 
As a $(0, m)$ tensor, $\nabla^{m} f$ can be viewed as a module map from 
the $C^{\infty}(U)$-module $\Gamma_{U}(T^{m}(TM^{\C}))$ to the 
$C^{\infty}(U)$-module $C^{\infty}(U)$. 
Since $\nabla^{m} f$ is linear in $f$, 
we can view $\nabla^{m}$ as a linear map from $C^{\infty}(U)$ to 
$\hom_{C^{\infty}(U))}(\Gamma_{U}(T^{m}(TM^{\C})), C^{\infty}(U))$. 
Since such a map corresponds to a linear map from 
$\Gamma_{U}(T^{m}(TM^{\C}))$ to $L(C^{\infty}(U))$,
we have a linear map
$$\psi^{m}_{U}: \Gamma_{U}(T^{m}(TM^{\C})) \to L(C^{\infty}(U))$$
corresponding to  $\left(\sqrt{-1}\right)^{m}\nabla^{m}$, 
where $L(C^{\infty}(U))$ is the space of all 
linear operators on $C^{\infty}(U)$. (Here we multiply 
$\nabla^{m}$ by the number $\left(\sqrt{-1}\right)^{m}$ 
because in canonical quantization, generalized momenta
should act on the space of functions by $\sqrt{-1}$ times differential 
operators.) By definition, 
for $\mathcal{X}\in \Gamma_{U}(T^{m}(TM^{\C}))$, 
$$(\psi^{m}_{U}(\mathcal{X}))f=\left(\sqrt{-1}\right)^{m}
(\nabla^{m}f)(\mathcal{X}).$$
The linear maps $\psi^{m}_{U}$ for $m\in \N$ 
give a single linear map 
$$\psi_{U}: \Gamma_{U}(T(TM^{\C})) \to L(C^{\infty}(U)).$$

We know that $\Gamma_{U}(T(TM^{\C}))$ is an 
associative algebra. The space $L(C^{\infty}(U))$ is in fact also an associative algebra.
But in general, the isomorphism $\psi_{U}$ is not an isomorphism of associative algebras. 
The associative algebra $\Gamma_{U}(T(TM^{\C}))$ has a subalgebra $\Pi_{U}(T(TM^{\C}))$. 
Let 
$$\phi_{U}: \Pi_{U}(T(TM^{\C})) \to L(C^{\infty}(U))$$
be the restriction of 
$\psi_{U}$ to $\Pi_{U}(T(TM^{\C}))$.
Then we have:

\begin{thm}\label{phi-hom}
For $\mathcal{X}\in \Gamma_{U}(T(TM^{\C}))$ and $\mathcal{Y}\in 
\Pi_{U}(T(TM^{\C}))$, we have
\begin{equation}\label{phi-hom-3}
\psi_{U}(\mathcal{X}\otimes \mathcal{Y})
=\psi_{U}(\mathcal{X})
\psi_{U}(\mathcal{Y}).
\end{equation}
In particular, the linear map $\phi_{U}$ 
is a homomorphism of associative algebras and
gives $C^{\infty}(U)$ a $\Pi_{U}(T(TM^{\C}))$-module structure. 
\end{thm}
\pf
We need only prove (\ref{phi-hom-3}) for
$m, l\in \N$,  $\mathcal{X}\in \Gamma_{U}(T^{m}(TM^{\C}))$ and $\mathcal{Y}\in 
\Pi_{U}(T^{m}(TM^{\C}))$. We use induction on $m$. 
When $m=0$, (\ref{phi-hom-3}) certainly holds. Now assume that when $m=k$, 
(\ref{phi-hom-3}) holds. To prove (\ref{phi-hom-3}) in the case $m=k+1$, 
we need only prove that for $f\in C^{\infty}(U)$ and $p\in U$,
\begin{equation}\label{phi-hom-1}
(\psi_{\tilde{U}}(\mathcal{X}\otimes \mathcal{Y})f)(p)
=(\psi_{U}(\mathcal{X})
\psi_{U}(\mathcal{Y})f)(p).
\end{equation}
For $p\in U$, there exists
an open subset $\tilde{U}$ of $U$ containing $p$ such that the restriction 
$\mathcal{X}|_{\tilde{U}}$ of $\mathcal{X}$ to $\tilde{U}$ is a sum of 
elements of the form $X\otimes \tilde{\mathcal{X}}$ for $X\in \Gamma_{\tilde{U}}(TM^{\C})$ and 
$\tilde{\mathcal{X}}\in \Pi_{\tilde{U}}(T^{k}(TM^{\C}))$. Hence we need
only prove (\ref{phi-hom-1})
for those $\mathcal{X}$ such that 
$$\mathcal{X}|_{\tilde{U}}=X\otimes \tilde{\mathcal{X}}$$
for 
$X\in \Gamma_{\tilde{U}}(TM^{\C})$ and 
$\tilde{\mathcal{X}}\in \Pi_{\tilde{U}}(T^{k}(TM^{\C}))$ 
In this case for $f\in C^{\infty}(U)$, 
by definition, 
\begin{eqnarray}\label{phi-hom-4}
\lefteqn{(\psi_{\tilde{U}}(\mathcal{X}|_{\tilde{U}})
\psi_{\tilde{U}}(\mathcal{Y}|_{\tilde{U}}))f}\nn
&&=(\psi_{\tilde{U}}(\mathcal{X}|_{\tilde{U}}))
((\psi_{\tilde{U}}(\mathcal{Y}|_{\tilde{U}}))f)\nn
&&=\left(\sqrt{-1}\right)^{k+1}
(\nabla^{k+1}((\psi_{\tilde{U}}(\mathcal{Y}|_{\tilde{U}}))f))
(X\otimes \tilde{\mathcal{X}})\nn
&&=\sqrt{-1}X((\left(\sqrt{-1}\right)^{k}\nabla^{k} 
((\psi_{\tilde{U}}(\mathcal{Y}|_{\tilde{U}}))f))
(\tilde{\mathcal{X}}))\nn
&&\quad -\sqrt{-1}(\left(\sqrt{-1}\right)^{k}\nabla^{k} 
((\psi_{\tilde{U}}(\mathcal{Y}|_{\tilde{U}})f))(\nabla_{X}\tilde{\mathcal{X}})\nn
&&=\sqrt{-1}X((\psi_{\tilde{U}}
(\tilde{\mathcal{X}}))
((\psi_{\tilde{U}}(\mathcal{Y}|_{\tilde{U}})f))
-\sqrt{-1}(\psi_{\tilde{U}}(\nabla_{X}\tilde{\mathcal{X}}))
((\psi_{\tilde{U}}(\mathcal{Y}|_{\tilde{U}}))f)\nn
&&=\sqrt{-1}X((\psi_{\tilde{U}}
(\tilde{\mathcal{X}})
\psi_{\tilde{U}}(\mathcal{Y}|_{\tilde{U}}))f)
-\sqrt{-1}(\psi_{\tilde{U}}(\nabla_{X}\tilde{\mathcal{X}})
\psi_{\tilde{U}}(\mathcal{Y}|_{\tilde{U}}))f.
\end{eqnarray}
By the induction assumption, the right-hand side of (\ref{phi-hom-4})
is equal to 
\begin{equation}\label{phi-hom-5}
\sqrt{-1}X((\psi_{\tilde{U}}
(\tilde{\mathcal{X}}\otimes \mathcal{Y}|_{\tilde{U}})f) -
\sqrt{-1}(\psi_{\tilde{U}}((\nabla_{X}\tilde{\mathcal{X}})
\otimes \mathcal{Y}|_{\tilde{U}})f.
\end{equation}
Since $\mathcal{Y}$  is parallel, we have 
$$\nabla_{X}\mathcal{Y}|_{\tilde{U}}=0$$
and thus
\begin{equation}\label{phi-hom-6}
(\nabla_{X}\tilde{\mathcal{X}})
\otimes \mathcal{Y}|_{\tilde{U}}=\nabla_{X}(\tilde{\mathcal{X}}
\otimes \mathcal{Y}|_{\tilde{U}}).
\end{equation}
Using (\ref{phi-hom-6}), (\ref{phi-hom-5}) becomes
\begin{eqnarray}\label{phi-hom-7}
\lefteqn{\sqrt{-1}X((\psi_{\tilde{U}}
(\tilde{\mathcal{X}}\otimes \mathcal{Y}|_{\tilde{U}})f) -
\sqrt{-1}(\psi_{\tilde{U}}(\nabla_{X}(\tilde{\mathcal{X}}
\otimes \mathcal{Y}|_{\tilde{U}})))f}\nn
&&=\sqrt{-1}X((\left(\sqrt{-1}\right)^{k+l}\nabla^{k+l} f)
(\tilde{\mathcal{X}}\otimes \mathcal{Y}|_{\tilde{U}}))\nn
&&\quad -
\sqrt{-1}(\left(\sqrt{-1}\right)^{k+l}\nabla^{k+l} f)
(\nabla_{X}(\tilde{\mathcal{X}}
\otimes \mathcal{Y)}|_{\tilde{U}})\nn
&&=(\left(\sqrt{-1}\right)^{k+1+l}\nabla^{k+1+l} f)
(X\otimes \tilde{\mathcal{X}}\otimes \mathcal{Y}|_{\tilde{U}})\nn
&&=(\psi_{\tilde{U}}
(\mathcal{X}|_{\tilde{U}}\otimes \mathcal{Y}|_{\tilde{U}}))f.
\end{eqnarray}

The calculations from (\ref{phi-hom-4}) 
to (\ref{phi-hom-7}) show that 
the left-hand side of (\ref{phi-hom-4}) is equal to the right-hand side of 
(\ref{phi-hom-7}). In particular, the value of the left-hand side of (\ref{phi-hom-4})
at $p$ is equal to the value of the right-hand side of 
(\ref{phi-hom-7}) at $p$. 
But the value of the left-hand side of (\ref{phi-hom-4}) at $p$ is equal to the 
right-hand side of (\ref{phi-hom-1}) and the value of the
right-hand side of (\ref{phi-hom-7}) at $p$ is equal to the left-hand side of 
(\ref{phi-hom-1}). Thus  (\ref{phi-hom-3}) holds. Since $p$ and $f$ are arbitrary,
(\ref{phi-hom-3}) in the case $m=k+1$ is proved. 
\epfv

\section{A sheaf $\mathcal{W}$ of modules 
for $\mathcal{V}$ constructed from the sheaf of smooth functions
on $M$}

In this section, we construct a sheaf $\mathcal{W}$ of left
modules for 
the sheaf $\mathcal{V}$ of meromorphic open-string 
vertex algebras from the sheaf 
of smooth functions on $M$.

Let $U$ be an open subset of $M$. For simplicity, we discuss only 
the case that $U$ is connected. The general case is similar. 
By Theorem \ref{phi-hom}, $C^{\infty}(U)$ 
is a $\Pi_{U}(T(TM^{\C}))$-module.
 For $p\in U$, by Proposition \ref{parallel-fixed-pts}, 
$\Pi_{U}(T(TM^{\C}))$ is isomorphic to 
$(T(T_{p}M^{\C}))^{H_{p}(U)}$. We shall identify  $\Pi_{U}(T(TM^{\C}))$ with
$(T(T_{p}M^{\C}))^{H_{p}(U)}$. In particular, 
$C^{\infty}(U)$ is a $(T(T_{p}M^{\C}))^{H_{p}(U)}$-module.
Since $(T(T_{p}M^{\C}))^{H_{p}(U)}$ is a subalgebra of $T(T_{p}M^{\C})$, 
we have the induced 
$T(T_{p}M^{\C})$-module  
$$C_{p}(U)=T(T_{p}M^{\C})\otimes_{(T(T_{p}M^{\C}))^{H_{p}(U)}}C^{\infty}(U).$$
By Theorems 5.1 (or Theorem \ref{buncle-mero-op-va} in Section 3 above)
and 6.5 in \cite{H3}, $T(\widehat{T_{p}M}_{-})$ has a natural structure
of meromorphic open-string vertex algebra and 
$T(\widehat{T_{p}M}_{-})\otimes C_{p}(U)$
has a natural structure of left $T(\widehat{T_{p}M}_{-})$-module. 
By Corollary \ref{holonomy}, $(T(\widehat{T_{p}M}_{-}))^{H_{p}(U)}$
is a meromorphic open-string vertex subalgebra of $T(\widehat{T_{p}M}_{-})$.
In particular, $T(\widehat{T_{p}M}_{-})\otimes C_{p}(U)$ is also a 
left $(T(\widehat{T_{p}M}_{-}))^{H_{p}(U)}$-module. Let 
$W^{0}_{U}$ be the left $(T(\widehat{T_{p}M}_{-}))^{H_{p}(U)}$-submodule of 
$T(\widehat{T_{p}M}_{-})\otimes C_{p}(U)$ generated by elements of the form
$1\otimes(1\otimes_{(T(T_{p}M^{\C}))^{H_{p}(U)}}f)$
for $f\in C^{\infty}(U)$, where $1\otimes_{(T(T_{p}M^{\C}))^{H_{p}(U)}}f$ 
is the image of $1\otimes f$ under the projection from 
$T(T_{p}M^{\C})\otimes C^{\infty}(U)$ to $C_{p}(U)$. 
By Theorem \ref{V-U-mero-op-va},
the meromorphic open-string vertex subalgebra 
$(T(\widehat{T_{p}M}_{-}))^{H_{p}(U)}$ is canonically isomorphic to 
$V_{U}=\Pi_{U}(T(\widehat{TM}_{-}))$. We shall identify 
$(T(\widehat{T_{p}M}_{-}))^{H_{p}(U)}$ and $V_{U}$.
Thus $W^{0}_{U}$ has 
a natural structure of left $V_{U}$-module.

The construction of $W^{0}_{U}$ here depends on $p$. But it is 
in fact independent of $p$. Let $q$ be another point in $U$.
Then the subspace of $T(\widehat{T_{p}M}_{-})\otimes C_{p}(U)$
consisting of elements of the form 
$1\otimes(1\otimes_{(T(T_{p}M^{\C}))^{H_{p}(U)}}f)$ 
for $f\in C^{\infty}(U)$ is canonically isomorphic to 
the subspace of $T(\widehat{T_{q}M}_{-})\otimes C_{q}(U)$
consisting of elements of the form 
$1\otimes(1\otimes_{(T(T_{q}M))^{H_{q}(U)}}f)$
for $f\in C^{\infty}(U)$.
Also both the meromorphic open-string vertex algebras
$(T(\widehat{T_{p}M}_{-}))^{H_{p}(U)}$ and $(T(\widehat{T_{q}M}_{-}))^{H_{q}(U)}$
are canonically isomorphic to 
the meromorphic open-string vertex algebra $\Pi_{U}(T(\widehat{TM}_{-}))=V_{U}$. 
Thus the  left $V_{U}$-module generated by elements of the form
$1\otimes(1\otimes_{(T(T_{p}M^{\C}))^{H_{p}(U)}}f)$
for $f\in C^{\infty}(U)$ is isomorphic to the left 
$V_{U}$-module generated by elements of the form
$1\otimes(1\otimes_{(T(T_{q}M))^{H_{q}(U)}}f)$
for $f\in C^{\infty}(U)$. Since the subspace of 
$T(\widehat{T_{p}M}_{-})\otimes C_{p}(U)$
consisting of elements of the form 
$1\otimes(1\otimes_{(T(T_{p}M^{\C}))^{H_{p}(U)}}f)$ 
for $f\in C^{\infty}(U)$ is canonically isomorphic to 
$C^{\infty}(U)$, we can view $W^{0}_{U}$ as a canonical left $V_{U}$-module 
generated by $C^{\infty}(U)$. 
We have proved the following result in the case that $U$ is 
connected; the general case
can be obtained using direct products as in the case of 
open-string vertex algebras
in Section 3:

\begin{thm}\label{W-U-left-module}
For an open subset $U$ of $M$ and $p\in U$, $W^{0}_{U}$ has 
a natural structure of left $V_{U}$-module. For different choices of $p$, 
we obtain canonically isomorphic left $V_{U}$-modules. In particular,
we have a canonical left $V_{U}$-module $W^{0}_{U}$ generated by 
$C^{\infty}(U)$ up to canonical isomorphisms. \epf
\end{thm}

Let $V_{1}$ and $V_{2}$ be meromorphic open-string vertex algebras  and 
$W_{1}$ and $W_{2}$ left $V_{1}$- and $V_{2}$-modules respectively. 
Let $f: V_{1}\to V_{2}$ be a homomorphism of meromorphic open-string vertex algebras.
Then $W_{2}$ is also a left $V_{1}$-module. 
A {\it homomorphism from $W_{1}$ to 
$W_{2}$ associated to $f$} is a homomorphism from $W_{1}$ to $W_{2}$ as 
a left $V_{1}$-module. 

By definition, $W^{0}_{U}$ is the left $(T(\widehat{T_{p}M}_{-}))^{H_{p}(U)}$-submodule of 
$T(\widehat{T_{p}M}_{-})\otimes C_{P}(U)$ generated by elements of the form
$1\otimes(1\otimes_{(T(T_{p}M^{\C}))^{H_{p}(U)}}f)$
for $f\in C^{\infty}(U)$. For open subsets $U$ and $\tilde{U}$ such that 
$\tilde{U}\subset U$, we have a restriction map from $V_{U}$ to $V_{\tilde{U}}$
which corresponding to the restriction map from $(T(\widehat{T_{p}M}_{-}))^{H_{p}(U)}$ to
$(T(\widehat{T_{p}M}_{-}))^{H_{p}(\tilde{U})}$. Using the restriction map 
from $(T(\widehat{T_{p}M}_{-}))^{H_{p}(U)}$ to
$(T(\widehat{T_{p}M}_{-}))^{H_{p}(\tilde{U})}$ and the restriction map from $C^{\infty}(U)$
to $C^{\infty}(\tilde{U})$, we obtain a restriction map from $W^{0}_{U}$ to 
$W^{0}_{\tilde{U}}$.
By definition, we obtain the following:

\begin{prop}
For open subsets $U$ and $\tilde{U}$ such that 
$\tilde{U}\subset U$, the restriction map from $W^{0}_{U}$ to 
$W^{0}_{\tilde{U}}$ is 
a homomorphism from $W^{0}_{U}$ to $W^{0}_{\tilde{U}}$ associated with the 
restriction map from $V_{U}$ to $V_{\tilde{U}}$. In particular, 
$U\to W^{0}_{U}$ gives a 
presheaf $\mathcal{W}^{0}$ of $\mathcal{V}$-modules. \epf
\end{prop}

Let $\mathcal{W}$ be the sheafification of the presheaf $\mathcal{W}^{0}$.
We denote the section of $\mathcal{W}$ on $U$ by $W_{U}$.
Then we obtain:

\begin{thm}\label{W}
The sheaf $\mathcal{W}$ is a sheaf of $\mathcal{V}$-modules.
In particular, the global section $W_{M}$ of $\mathcal{W}$ on $M$ 
is a left $V_{M}$-module such that $C^{\infty}(M)$ is embedded 
as a subspace of $W_{M}$.
\end{thm}
\pf 
We prove only the statement that $C^{\infty}(M)$ is embedded 
as a subspace of $W_{M}$. As we have mentioned above, 
for an open subset $U$ of $M$, 
$f\mapsto 1\otimes(1\otimes_{(T(T_{p}M^{\C}))^{H_{p}(U)}}f)$
for $f\in C^{\infty}(U)$ give an injective linear map from 
$C^{\infty}(U)$ to $W^{0}_{U}$. Thus we can identify 
$f$ and $1\otimes(1\otimes_{(T(T_{p}M^{\C}))^{H_{p}(U)}}f)$
and view $C^{\infty}(U)$
as a subspace of $W^{0}_{U}$. After this identification, we see that
$W^{0}_{U}$ is generated by $C^{\infty}(U)$. In particular, 
$W^{0}_{M}$ is generated by $C^{\infty}(M)$. 
From the construction of the sheafification of a
presheaf, we see that for an open subset $U$ of $M$, 
$W^{0}_{U}$ can be embedded into $W_{U}$. Thus $C^{\infty}(U)$ can also
be viewed as a subspace of $W_{U}$. In particular, $C^{\infty}(M)$
can be viewed as a subspace of $W_{M}$.
\epfv

\section{The Laplacian on $M$ as a component of a vertex operator}

Many conjectures in geometry were obtained from quantum field theory by 
interpreting some geometric or analytic objects as quantum-field-theoretic 
objects. In this section, as an example, we show that the Laplacian of 
$M$ is in fact a component of a vertex operator for $W_{M}$ acting on $C^{\infty}(M)$. 

Let $\{E_{i}\}_{i=1}^{n}$ be an 
orthonormal frame in an open neighborhood $U$ of a point $p\in M$. 
Recall the element $g_{\C}^{-1}(-1, -1)\in V_{U}$. Then in $U$, 
$$g_{\C}^{-1}(-1, -1)=\sum_{i=1}^{n}(E_{i}\otimes_{\R} t^{-1})\otimes 
(E_{i}\otimes_{\R} t^{-1}).$$
We identify $V_{U}=\Pi_{U}(T(\widehat{TM}_{-})$ with $(T(\widehat{T_{p}M}_{-}))^{H_{p}(U)}$.
Under this identification, $g_{\C}^{-1}(-1, -1)$ is identified with 
$$\sum_{i=1}^{n}(E_{i}|_{p}\otimes_{\R} t^{-1})
\otimes (E_{i}|_{p}\otimes_{\R} t^{-1})
=\sum_{i=1}^{n}(E_{i}|_{p})(-1)(E_{i}|_{p})(-1)\one,$$
where $\one_{T(\widehat{T_{p}M}_{-})}$
is the vacuum of the meromorphic open-string vertex algebra 
$T(\widehat{T_{p}M}_{-})$ and $(E_{i}|_{p})(-1)$ is the representation image of 
$E_{i}|_{p}\otimes_{\R} t^{-1}$ on $T(\widehat{T_{p}M}_{-})$. 

Recall that $W^{0}_{U}$  by definition is 
the left $(T(\widehat{T_{p}M}_{-}))^{H_{p}(U)}$-submodule of 
$T(\widehat{T_{p}M}_{-})\otimes C_{p}(U)$ generated by elements of the form
$1\otimes(1\otimes_{(T(T_{p}M^{\C}))^{H_{p}(U)}}f)$
for $f\in C^{\infty}(U)$. Let $f\in C^{\infty}(U)$. Consider 
$$Y_{W^{0}_{U}}(-g_{\C}^{-1}(-1, -1), x)
(1\otimes(1\otimes_{(T(T_{p}M^{\C}))^{H_{p}(U)}}f)).$$
Since  $W^{0}_{U}$ is a $(T(\widehat{T_{p}M}_{-}))^{H_{p}(U)}$-submodule of 
$T(\widehat{T_{p}M}_{-})\otimes C_{p}(U)$ and the 
$(T(\widehat{T_{p}M}_{-}))^{H_{p}(U)}$-module
$T(\widehat{T_{p}M}_{-})\otimes C_{p}(U)$ is induced from the 
$T(\widehat{T_{p}M}_{-})$-module structure on $T(\widehat{T_{p}M}_{-})\otimes C_{p}(U)$,
the vertex operator map $Y_{W^{0}_{U}}$ is the restriction of the vertex operator map
$Y_{T(\widehat{T_{p}M}_{-})}$ to $(T(\widehat{T_{p}M}_{-}))^{H_{p}(U)}\otimes 
W^{0}_{U}$. In particular, 
\begin{eqnarray*}
Y_{W^{0}_{U}}(-g_{\C}^{-1}(-1, -1), x)
&=&Y_{T(\widehat{T_{p}M}_{-})}(-g_{\C}^{-1}(-1, -1), x)\nn
&=&-\sum_{i=1}^{n} Y_{T(\widehat{T_{p}M}_{-})}((E_{i}|_{p})(-1)(E_{i}|_{p})(-1)\one, x)\nn
&=&-\sum_{i=1}^{n} \nord (\sqrt{-1}E_{i}|_{p})(x) (\sqrt{-1}E_{i}|_{p})(x)\nord\nn
&=&\sum_{i=1}^{n} \nord (E_{i}|_{p})(x) (E_{i}|_{p})(x)\nord.
\end{eqnarray*}
Then the coefficient of the $x^{-2}$ term of $Y_{W_{U}}(-g_{\C}^{-1}(-1, -1), x)$
is 
\begin{eqnarray*}
\lefteqn{\sum_{i=1}^{n} \sum_{k\in \Z} \nord (E_{i}|_{p})(-k) (E_{i}|_{p})(k)\nord}\nn
&&=2\sum_{i=1}^{n} \sum_{k\in \Z_{+}} (E_{i}|_{p})(-k) (E_{i}|_{p})(k)
+\sum_{i=1}^{n} (E_{i}|_{p})(0) (E_{i}|_{p})(0).
\end{eqnarray*}
This coefficient or component of the vertex operator $Y_{W_{U}}(g_{\C}^{-1}(-1, -1), x)$
acting on $(1\otimes(1\otimes_{(T(T_{p}M^{\C}))^{H_{p}(U)}}f))$
is equal to 
\begin{eqnarray*}
\lefteqn{\sum_{i=1}^{n} (E_{i}|_{p})(0) (E_{i}|_{p})(0)
(1\otimes(1\otimes_{(T(T_{p}M^{\C}))^{H_{p}(U)}}f))}\nn
&&=1\otimes\left(\sum_{i=1}^{n}(E_{i}|_{p}\otimes E_{i}|_{p}) 
\otimes_{(T(T_{p}M^{\C}))^{H_{p}(M)}}f\right)\nn
&&=1\otimes\left(1\otimes_{(T(T_{p}M^{\C}))^{H_{p}(U)}}
\left(\phi\left(\sum_{i=1}^{n}E_{i}\otimes E_{i}\right) f\right)\right)\nn
&&=1\otimes\left(1\otimes_{(T(T_{p}M^{\C}))^{H_{p}(U)}}
\left((\nabla^{2}f)\left(\sum_{i=1}^{n}E_{i}\otimes E_{i}\right)\right)\right)\nn
&&=1\otimes\left(1\otimes_{(T(T_{p}M^{\C}))^{H_{p}(U)}}
\left(\sum_{i=1}^{n}\nabla^{2}_{E_{i}, E_{i}}f\right)\right)\nn
&&=1\otimes (1\otimes_{(T(T_{p}M^{\C}))^{H_{p}(U)}}
\Delta f).
\end{eqnarray*}
If we identify $f$ with $1\otimes(1\otimes_{(T(T_{p}M^{\C}))^{H_{p}(U)}}f)$,
then we see that this component of the vertex operator 
$Y_{W^{0}_{U}}(-g_{\C}^{-1}(-1, -1), x)$
acting on $f$ is equal to 
the Laplacian $\Delta$ acting on $f$. 

Since $\mathcal{V}$ and 
$\mathcal{W}$ are sheaves and $\mathcal{W}$ is the sheafification of the 
presheaf $\mathcal{W}^{0}$, the conclusion 
above 
for $\mathcal{W}^{0}$ and smooth functions 
on a small open neighborhood of 
every $p\in M$ implies that the same conclusion holds for 
$\mathcal{W}$ and smooth functions 
on any open subset $U$ of $M$. In particular, if we view 
$f\in C^{\infty}(M)$ as an element of $W_{M}$, then the 
component of the vertex operator $Y_{W_{M}}(-g_{\C}^{-1}(-1, -1), x)$
acting on $f$
is equal to $\Delta f$.

\noindent {\small \sc Department of Mathematics, Rutgers University,
110 Frelinghuysen Rd., Piscataway, NJ 08854-8019}

\vspace{1em}
\noindent {\it and}

\vspace{1em}
\noindent {\small \sc Beijing International Center for Mathematical Research,
Peking University, Beijing, China}
\vspace{1em}

\noindent {\em E-mail address}: yzhuang@math.rutgers.edu


\begin{thebibliography}{HHH}


\bibitem[H1]{H1}
Y.-Z. Huang, {\em Two-dimensional Conformal Geometry and Vertex
Operator Algebras}, Progress in Math., Vol. 148, Birkh\"{a}user,
Boston, 1997.

\bibitem[H2]{H2}
Y.-Z. Huang, Intertwining operator algebras, genus-zero modular functors 
and genus-zero conformal field theories, in: 
{\it Operads: Proceedings of Renaissance Conferences}, 
ed. J.-L. Loday, J. Stasheff, and A. A. Voronov, 
Contemporary Math., Vol. 202, Amer. Math. Soc., Providence, 1997, 335--355. 

\bibitem[H3]{H3}
Y.-Z. Huang, Meromorphic open-string vertex algebras, {\it J. Math. Phys.} {\bf 54} (2013), 051702. 

\bibitem[HK]{HK} Y.-Z. Huang and L. Kong,
Open-string vertex algebras, tensor categories and operads, 
{\it Comm. Math. Phys.} {\bf  250} (2004), 433--471. 


\bibitem[P]{P}
P. Petersen, {\it Riemannian Geometry}, Graduate Texts in Math., Vol. 171, Springer,
New York, 1997.

\bibitem[Q1]{Q1}
F. Qi, Meromorphic open-string vertex algebras and modules 
over two-dimensional orientable space forms,
{\it Lett. Math. Phys.} {\bf 111} (2021), article no. 27. 

\bibitem[Q2]{Q2}
F. Qi, Covariant derivatives of eigenfunctions along parallel
tensors over space forms and a conjecture motivated by
the vertex algebraic structure,
{\it J. Noncommut. Geom.} {\bf 16} (2022), 717--759.

\bibitem[T1]{T1}
H. Tamanoi, {\it Elliptic genera and vertex operator super algebras}, {\it Proc.
Japan Acad.} {\bf A71} (1995),  177--181.

\bibitem[T2]{T}
H. Tamanoi, {\it Elliptic genera and vertex operator super-algebras},
Lecture Notes in Mathematics, Vol. 1704, Springer-Verlag, Berlin, 1999. 


\end{thebibliography}
\end{document}